\documentclass[a4,11pt,twoside,leqno]{article}
\usepackage{amsfonts}
\usepackage{amssymb}
\usepackage{amsthm}
\usepackage{amsmath}

\newtheorem{thm}{Theorem}[section]
\newtheorem{cor}[thm]{Corollary}
\newtheorem{lem}[thm]{Lemma}

\theoremstyle{definition}
\newtheorem{defn}[thm]{Definition}
\theoremstyle{remark}
\newtheorem{rem}[thm]{Remark}

\numberwithin{equation}{section}

\begin{document}

\title{Harnack inequalities  for graphs with non-negative Ricci curvature}

\author{Fan Chung\\University of California, San Diego\
\and Yong Lin \\ Renmin University of China \and S.-T. Yau
\\Harvard University }

\maketitle

\footnotetext{Classification:[2000]{31C20, 31C05}\\

Keywords: The Laplace operator for graphs, the Harnack inequalities, eigenvalues, diameter.}

\vskip0.3cm

\begin{abstract}
We establish a Harnack inequality for finite connected graphs with non-negative Ricci
curvature. As a consequence, we derive an eigenvalue lower bound,  extending previous results for
 Ricci flat graphs.  
\end{abstract}

\section{Introduction}

 Let $G$ be an undirected finite
connected weighted graph with vertex set $V$ and edge set $E$.  The edge weight of an edge $\{x,y\}$ is denoted by
$w_{xy}$ and the degree $d_x$ is the sum of all $w_{xy}$ over all $y$ adjacent to $x$.
The Laplace operator $\Delta$ of a graph $G$ is defined by
$$\Delta{f(x)}=\frac{1}{d_x}\sum_{y\sim{x}}w_{xy}\big(f(y)-f(x)\big)$$
for  any function $f\in{V^R}=\{f|f:V\rightarrow{R}\}$ and any vertex $x \in V$ .

Suppose a function $f\in{V^R}$ satisfies that, for every vertex $x \in V$,
$$(-\Delta)f(x)=\frac{1}{d_x}\sum_{y\sim{x}}w_{xy}\big(f(x)-f(y)\big)=\lambda{f(x)}.$$
Then $f$ is called a harmonic eigenfunction of the Laplace operator
$\Delta$ on $G$ with eigenvalue $\lambda$. For a finite graph, it is straightforward to verify that $\lambda$  is
an eigenvalue for the (normalized) Laplacian  $\mathcal L$ as a matrix defined by
$${\mathcal L}= - \Delta = 1-D^{-1/2}AD^{-1/2} $$
where $D$ is the diagonal degree matrix and $A$ is the weighted adjacency matrix with $A(x,y)=w_{xy}$. Because of the positivity of $\mathcal L$, a connected graph has all eigenvalues positive except for one eigenvalue $zero$.

From this definition, at each vertex $x$, the eigenfunction locally stretches the
incident edges in a balanced fashion. Globally, it is desirable to have
some tools to capture the notion  that adjacent vertices are close to each other.

A crucial part of spectral graph theory concerns understanding the
behavior of eigenfunctions. Harnack inequalities are one of the main methods
for dealing with eigenfunctions. The  Harnack
inequalities  for certain special families of graphs, called Ricci flat graphs (see \cite{CHUNG-1994}), are
formulated as follows:\\
\begin{equation}
\frac{1}{d_x}\sum_{y\sim{x}}(f(x)-f(y))^2\leq{8}\lambda\max_zf^2(z)
\end{equation}
for   any
eigenfunction $f$ with eigenvalue $\lambda > 0$.

In general, the above inequality does not hold for all graphs.
An easy counterexample is the graph formed by joining two complete
graphs of the same size by a single edge  \cite{CHUNG-1994}.

In this paper, we will establish  a Harnack  inequality for
general graphs. We consider graphs with non-negative Ricci curvature $\kappa$ and we will show the following:
\begin{equation}
\frac{1}{d_x}\sum_{y\sim{x}}(f(x)-f(y))^2\leq{(8 \lambda -4\kappa)}\max_zf^2(z)
\end{equation}
for any graph with Ricci curvature $\kappa$. The definition of the Ricci curvature for graphs will be given in the next section.

For a graph $G$, the diameter of $G$ is the least number $D$ such that any two vertices in $G$ are joined by a path with at most $D$ edges.
By using the above Harnack inequality, we will derive
 the following eigenvalue/diameter inequality  for graphs with non-negative
Ricci curvature.
$$\lambda\geq\frac{1+2\kappa dD^2}{4d\cdot{D}^2}$$
where $d$ is the maximum degree and $D$ denotes the diameter of $G$ (i.e., any two vertices can be joined by a path with at most $D$ edges).

\section{The Ricci curvature for graphs}
\indent In   \cite{CHUNG-1994} and \cite{CHUNG-1996}, Chung
and Yau defined Ricci flat graphs and proved that
inequality (1.1) and (1.4) hold for
a large family of Ricci flat graphs.
There are several ways to define Ricci curvature for a general graph.
In this paper, we will use the definition of
 Ricci curvature for graphs
in the sense of Bakry-Emery \cite{BE-1985}, as introduced in \cite{LY}.  We note that a different notion of
Ricci curvature was introduced by  Ollivier \cite{olli-2007}.

To define the Ricci curvature of a graph, we begin with  a bilinear operator $\Gamma: V^R\times V^R\to
V^R,$ defined by
$$\Gamma(f,g)(x)=\frac 12\{\Delta(f(x)g(x))-f(x)\Delta g(x)-g(x)\Delta f(x)\}.$$ According to Bakry and
Emery \cite{BE-1985}, the Ricci curvature operator $\Gamma_2$ is
defined by:
$$\Gamma_2(f,g)(x)=\frac 12\{\Delta\Gamma(f,g)(x)-\Gamma(f,\Delta g)(x)-\Gamma(g,\Delta f)(x)\}.$$

For simplicity, we will omit the variable $x$ in the following
equations. Note that all the equations  hold locally
for every $x\in V$.

\begin{defn}
The operator $\Delta$ satisfies the curvature-dimension type
inequality $CD(m,\kappa) $ for $m\in (1,+\infty) $ if $$\Gamma_2(f,f)\ge \frac
1m(\Delta f)^2+\kappa\Gamma (f,f).$$
\end{defn}
We call $m$ the dimension of the operator $\Delta$ and $\kappa$ a lower
bound of the Ricci curvature of the operator $\Delta$. If $\Gamma _2\ge \kappa\Gamma$, we say that $\Delta$ satisfies $CD(\infty
,\kappa)$.

It is easy to see that for $m<m^{\prime}$, the operator $\Delta$
satisfies the curvature-dimension type inequality $CD(m^{\prime},K)$
if it satisfies the curvature-dimension type inequality $CD(m,K)$.

Here we list a number of helpful facts concerning $\Gamma$, $\Gamma_2$ and the Ricci curvature that will be
useful later.

From the definition of $\Gamma$, we can  express $\Gamma$ in the following alternative formulation. The derivation is straightforward and we omit the proof here.

\begin{lem}
\label{m1}
\begin{eqnarray}
\Gamma (f,g)(x)&=&{1\over {2d_x}}\sum_{y\sim{x}}w_{xy}\big(f(x)-f(y)\big) \big(g(x)-g(y)\big).\label{eq3}\\
\Gamma (f,f)(x)&=&{1\over 2d_x}\sum_{y\sim{x}}w_{xy}[f(x)-f(y)]^2={1\over 2}\vert \nabla
f\vert^{2}(x). \label{eq4}
\end{eqnarray}
\end{lem}

For Laplace-Beltrami operator $\Delta$ on a complete $m$ dimensional
Riemannion manifold, the operator $\Delta$ satisfies $CD(m,K)$ if 
the Ricci curvature of the Riemanian manifold is bounded below by
 $\kappa$. For graphs, a similar bound can be established as follows.
\begin{lem}
\label{m2}
In a connected graph $G$, let $\lambda$ denote a non-trivial eigenvalue. Then the Ricci curvature $\kappa$ of $G$
satisfies
$$\kappa \leq \lambda.$$
\end{lem}
Proof: Let $f$ denote a harmonic eigenvector associated with eigenvalue $\lambda$.
Consider the vertex $x$ which achieves the maximum of $\vert \nabla
f\vert^{2}$. Then we have $\Delta \vert \nabla
f\vert^{2} \leq 0$ and therefore
\begin{eqnarray*}
\Gamma_2(f,f) &\leq& -\Gamma(f, \Delta f)\\
&=& \lambda \Gamma(f,f).
\end{eqnarray*}
From the definition of $\kappa$, we have
$$  \lambda \Gamma(f,f) \geq \frac
1m(\Delta f)^2+\kappa\Gamma (f,f).$$
Thus we have $\kappa \leq \lambda$.
\hfill $\square$

We remark that Lemma \ref{m2} can be slightly improved to $\lambda \geq \kappa (1+ 1/(m-1))$ as seen in \cite{cll10}.

It was proved in \cite{LY} that the Ricci flat graphs as defined in
\cite{CHUNG-1994} and \cite{CHUNG-1996} are graphs satisfy $CD(\infty,0)$.
In \cite{LY}, it was  shown that any locally finite
connected $G$ satisfy the $CD(\frac{1}{2},\frac{1}{d}-1)$, if  the maximum degree $d$ is
finite, or $CD(2,-1)$ if $d$ is infinite. Thus, the Ricci curvature
of  a graph $G$ has a lower bound $-1$.



\section{Harnakc inequality and eigenvalue estimate}

First, we will establish several basic facts  for graphs with non-negative Ricci curvature.\\

\begin{lem}\label{lem1}
Suppose $G$ is a   finite connected graph satisfying $CD(m,\kappa)$. Then
for $x\in{V}$ and $f \in V^{R}$, we have
$$(\frac{4}{m}-2)(\Delta{f}(x))^2+(2+2\kappa)|\nabla{f}|^2(x)\leq\frac{1}{d_x}
\sum_{y\sim{x}}\frac{w_{xy}}{d_y}\sum_{z\sim{y}}w_{yz}[f(x)-2f(y)+f(z)]^2.$$
\end{lem}
Proof:  We consider $\Delta(\Gamma(f,f))$. By straightforward manipulation and (\ref{eq4}), we have

\begin{eqnarray*}
\Delta(\Gamma(f,f))(x)&=&\frac{1}{2d_x}\sum_{y\sim{x}}\frac{w_{xy}}{d_y}\sum_{z\sim{y}}w_{yz}\big(-[f(x)-f(y)]^2+[f(y)-f(z)]^2\big)\\
&=&\frac{1}{2d_x}\sum_{y\sim{x}}\frac{w_{xy}}{d_y}\sum_{z\sim{y}}w_{yz}[f(x)-2f(y)+f(z)]^2\\
& &
-\frac{1}{d_x}\sum_{y\sim{x}}\frac{w_{xy}}{d_y}\sum_{z\sim{y}}w_{yz}[f(x)-2f(y)+f(z)][f(x)-f(y)],
\end{eqnarray*}
and by (\ref{eq3}) we have
\begin{eqnarray*}
\Gamma(f,\Delta f)(x)=\frac 12 \cdot {1\over d_x}\sum_{y\sim{x}}w_{xy}[f(y)-f(x)]\cdot [\Delta f(y)-\Delta f(x)].
\end{eqnarray*}

 By the definition of $\Gamma_2(f,f)$, we have

\begin{eqnarray}
\Gamma_2(f,f)(x)&=&\frac{1}{4}\frac{1}{d_x}\sum_{y\sim{x}}\frac{w_{xy}}{d_y}\sum_{z\sim{y}}w_{yz}[f(x)-2f(y)+f(z)]^2\nonumber\\
& &
-\frac{1}{2}\frac{1}{d_x}\sum_{y\sim{x}}w_{xy}[f(x)-f(y)]^2+\frac{1}{2}\left(\frac{1}
{d_x}\sum_{y\sim{x}}w_{xy}(f(x)-f(y))\right)^2\nonumber\\
&=&\frac{1}{4}\frac{1}{d_x}\sum_{y\sim{x}}\frac{w_{xy}}{d_y}\sum_{z\sim{y}}w_{yz}[f(x)-2f(y)+f(z)]^2- frac 1 2 \vert \nabla f \vert^2(x) + \frac 1 2 (\Delta f)^2. \label{eq5}
\end{eqnarray}

Since $G$ satisfies  $CD(m,\kappa)$, we have

$$\Gamma_2(f,f)\geq\frac{1}{m}(\Delta{f})^2+\kappa \Gamma(f,f).$$

From above inequality, we obtain

$$(\frac{1}{m}-\frac{1}{2})(\Delta{f})^2+\frac{1+\kappa}{2}\vert \nabla
f\vert^{2} \leq\frac{1}{4}
\frac{1}{d_x}\sum_{y\sim{x}}\frac{w_{xy}}{d_y}
\sum_{z\sim{y}}w_{yz}[f(x)-2f(y)+f(z)]^2$$
as desired.
\hfill $\square$

By using Lemma \ref{lem1}, we can prove the following Harnack type
inequality. The idea of proof comes from \cite{CHUNG-1994}.

\begin{thm}\label{thm1}
Suppose that a finite connected graph $G$ satisfies \ $CD(m,\kappa)$ and
$f\in{V^R}$ is a harmonic eigenfunction of Laplacian $\Delta$ with
eigenvalue $\lambda$. Then the following inequality holds for
all $x\in{V}$ and $\alpha\geq 2-2\kappa/\lambda$
$$|\nabla{f}|^2(x)+\alpha\lambda{f}^2(x)\leq\frac{({\alpha}^2-\frac{4}{m})
\lambda+2\kappa \alpha}{(\alpha-2)\lambda+2\kappa }\lambda\max_{z\in{V}}f^2(z).$$
\end{thm}

Proof:
Using Lemmas \ref{m1} and  \ref{lem1}, we have

\begin{eqnarray*}
(-\Delta)|\nabla{f}|^2(x)&=&-\frac{1}{d_x}\sum_{y\sim{x}}\frac{w_{xy}}{d_y}\sum_{z\sim{y}}w_{yz}[f(x)-2f(y)+f(z)]^2\\
& &
+\frac{2}{d_x}\sum_{y\sim{x}}\frac{w_{xy}}{d_y}\sum_{z\sim{y}}w_{yz}[f(x)-2f(y)+f(z)]\cdot[f(x)-f(y)]\\
&\le&-(2+2\kappa)\cdot|\nabla{f}|^2(x)+(2-\frac{4}{m})\cdot(\Delta{f}(x))^2+2\cdot|\nabla{f}|^2(x)\\
& &
+\frac{2}{d_x}\sum_{y\sim{x}}w_{xy}(f(x)-f(y))\cdot\frac{1}{d_y}\sum_{z\sim{y}}w_{yz}[f(z)-f(y)]\\
&=&-2\kappa \cdot|\nabla{f}|^2(x)+(2-\frac{4}{m})[-\lambda{f}(x)]^2+\frac{2}{d_x}\sum_{y\sim{x}}w_{xy}(f(x)-f(y))\cdot(-\lambda{f}(y))\\
&=&-2\kappa\cdot|\nabla{f}|^2(x)+(2-\frac{4}{m}){\lambda}^2f^2(x)\\
& &
+\frac{2}{d_x}\sum_{y\sim{x}}w_{xy}(f(x)-f(y))(-\lambda{f}(y)+\lambda{f}(x)-\lambda{f}(x))\\
&=&(2\lambda-2\kappa)\cdot|\nabla{f}|^2(x)-\frac{4}{m}{\lambda}^2f^2(x).
\end{eqnarray*}
Now we consider
\begin{eqnarray*}
(-\Delta)f^2(x)&=&\frac{1}{d_x}\sum_{y\sim{x}}w_{xy}[f^2(x)-f^2(y)]\\
&=&\frac{2}{d_x}\sum_{y\sim{x}}w_{xy} {f}(x)[f(x)-f(y)]-\frac{1}{d_x}\sum_{y\sim{x}}w_{xy}[f(x)-f(y)]^2\\
&=&2\lambda{f}^2(x)-|\nabla{f}|^2(x).
\end{eqnarray*}

Combining the above inequalities,  we have, for any positive $\alpha$, the following:

\begin{eqnarray*}
(-\Delta)(|\nabla{f}|^2(x)+\alpha\lambda{f^2}(x))&\leq&(2\lambda-2\kappa)|\nabla{f}|^2(x)-\frac{4}{m}
{\lambda}^2f^2(x)+2\alpha{\lambda}^2f^2(x)-\alpha\lambda|\nabla{f}|^2(x)\\
&=&(2\lambda-\alpha\lambda-2\kappa)|\nabla{f}|^2(x)+(2\alpha-\frac{4}{m}){\lambda}^2f^2(x).
\end{eqnarray*}

We choose a vertex $v$, which maximizes the expression
$$|\nabla{f}|^2(x)+\alpha\lambda{f}^2(x)$$
over all $x\in{V}.$
Then we have

\begin{eqnarray*}
0&\leq&(-\Delta)(|\nabla{f}|^2(v)+\alpha\lambda{f^2}(v))\\
&\leq&(2\lambda-\alpha\lambda-2\kappa)\cdot|\nabla{f}|^2(v)+(2\alpha-\frac{4}{m}){\lambda}^2f^2(v).
\end{eqnarray*}

This implies

$$|\nabla{f}|^2(v)\leq\frac{2\alpha-\frac{4}{m}}{(\alpha-2)\lambda+2\kappa}\cdot{\lambda}^2\cdot{f^2}(v)$$

for $\alpha>2-\frac{2\kappa}{\lambda}$.\\

Therefore for every $x\in{V}$, we have

\begin{eqnarray*}
|\nabla{f}^2|(x)+\alpha\lambda{f^2}(x)&\leq&|\nabla{f}|^2(v)+\alpha\lambda{f^2(v)}\\
&\leq&\frac{2\alpha-\frac{4}{m}}{(\alpha-2)\lambda+2\kappa}\cdot{\lambda}^2f^2(v)+\alpha\lambda{f}^2(v)\\
&\leq&\frac{({\alpha}^2-\frac{4}{m})\lambda+2\kappa \alpha}{(\alpha-2)\lambda+2\kappa}\cdot\lambda\cdot\max_{z\in{V}}{f}^2(z)
\end{eqnarray*}
as desired.
\hfill $\square$

From Lemma \ref{m2}, we can choose $\alpha=4-\frac{2\kappa}{\lambda}\geq 0$. By substituting into the statement of
Theorem \ref{thm1}, we have

\begin{thm}\label{thm2}
Suppose a finite connected graph $G$ satisfies the $CD(m,\kappa)$ and
$f\in{V^R}$ is a harmonic eigenfunction of Laplacian $\Delta$ with
nontrivial eigenvalue $\lambda$. Then the following inequality holds for all
$x\in{V}$
$$|\nabla{f}|^2(x)\leq \big((8-\frac{2}{m})\lambda-4\kappa \big)
\cdot\max_{z\in{V}}f^2(z).$$
\end{thm}
If $G$ is a non-negative Ricci curvature graph, i.e. $\kappa=0$. Then we
have the following result:

\begin{cor}\label{cor1}
Suppose a finite connected graph $G$ satisfies $CD(m,\kappa)$ and
$f\in{V^R}$ is a harmonic eigenfunction of Laplacian $\Delta$ with
nontrivial eigenvalue $\lambda$. Then the following Harnack inequality holds
for all $x\in{V}$
$$|\nabla{f}|^2(x)\leq(8-\frac{2}{m})\cdot\lambda\cdot\max_{z\in{V}}f^{2}(z).$$
\end{cor}
We can use the Harnack inequality in Theorem \ref{thm2} to derive
the following eigenvalue estimate.

\begin{thm}\label{thm3}
Suppose a finite connected graph $G$ satisfies $CD(m,\kappa)$ and  $\lambda$
is a non-zero eigenvalue of Laplace operator $\Delta$ on $G$. Then
$$\lambda\geq\frac{1+4\kappa dD^2}{d\cdot (8-\frac{2}{m})\cdot{D}^2}$$
where $d$ is the maximum degree and $D$ denotes the diameter of $G$.
\end{thm}
Proof: Let $f$ be the eigenfunction of Laplacian $\Delta$ with
eigenvalue $\lambda\neq{0}$. That is, for all $x \in V$,
$$(-\Delta)f(x)=\lambda{f}(x),$$  Then
\begin{eqnarray*}
\sum_{x\in{V}}d_xf(x)&=&\frac{1}{\lambda}\sum_{x\in{V}}d_x(-\Delta)f(x)\\
&=&\frac{1}{\lambda}\sum_{x\in{V}}\sum_{y\sim{x}}w_{xy}[f(x)-f(y)]\\
&=&0.
\end{eqnarray*}
We can assume that
$$\sup_{z\in{V}}f(z)=1>\inf_{z\in{V}}f(z)=\beta,$$
where $\beta<0$.\\
\indent Choose $x_1,x_t\in{V}$ such that
$f(x_1)=\sup_{z\in{V}}f(z)=1$, $f(x_t)=\inf_{z\in{V}}f(z)=\beta<0$ and
let $x_1$, $x_2$,...,$x_t$ be the shortest path connecting $x_1$ and
$x_t$, where $x_i\sim x_{i+1}$. Then $n\leq{D}$ where $D$ is the diameter
of $G$. From the Corrolalry \ref{cor1} we have

$$[f(x_{i-1})-f(x_{i})]^2+[f(x_i)-f(x_{i+1})]^2\leq d\cdot |\nabla{f}|^2(x_i)\leq d\cdot \big((8-\frac{2}{m})\cdot\lambda-4\kappa\big).$$

Therefore
$$\sum_{i=0}^{t-1}[f(x_i)-f(x_{i+i})]^2\leq dD\cdot\big((4-\frac{1}{m})\lambda-2\kappa \big).$$
On the other hand, by using the Cauchy-Schwarz inequality we have

\begin{eqnarray*}
\sum_{i=0}^{t-1}[f(x_i)-f(x_{i+1})]^2&\geq&\frac{1}{D}(f(x_t)-f(x_1))^2\\
&\geq&\frac{1}{D}
\end{eqnarray*}

Together we have

$$\lambda\geq\frac{1+2\kappa dD^2}{d\cdot (4-\frac{1}{m})\cdot{D}^2}.$$

This completes the proof of Theorem \ref{thm3}.
\hfill $\square$

\begin{rem}

We note that  Theorem \ref{thm3} gives an eigenvalue lower bound for graphs with
non-negative Ricci curvature $\kappa $ satisfying
$$ \kappa >-\frac{1}{2dD^2}.$$

\end{rem}

As an immediately consequence, we have the following:
\begin{cor}\label{cor3}
Suppose a finite connected graph $G$ satisfies $CD(m,0)$ and  $\lambda$
is a non-zero eigenvalue of Laplace operator $\Delta$ on $G$. Then
$$\lambda\geq\frac{1}{d\cdot (4-\frac{}{m})\cdot{D}^2}$$
where $d$ is the maximum degree and $D$ denotes the diameter of $G$.
\end{cor}

Chung and Yau(see \cite{CHUNG-1994}) proved that $\lambda\geq\frac{1}{d\cdot 8\cdot{D}^2}$ for
Ricci flat graphs. Since Ricci flat graphs satisfies $CD(\infty,0)$,
 our results extend  and strengthen the  results in  \cite{CHUNG-1994}.

\end{document}